\documentclass[a4paper,12pt,reqno]{amsart}
\usepackage{amsfonts}
\usepackage[colorlinks,linkcolor=blue,citecolor=blue]{hyperref}
\usepackage{latexsym, amssymb, amsmath, amsthm, bbm}
\usepackage[all]{xy}
\usepackage{enumerate}
\usepackage{mathrsfs}
\usepackage{geometry}
\usepackage{pgfplots}

\def \k{\Bbbk}

\usepackage{mathtools}

\numberwithin{equation}{section}
\theoremstyle{plain}
 \newtheorem{theorem}{Theorem}[section]

\newtheorem{lemma}[theorem]{Lemma}

\pgfplotsset{compat=1.18}

\begin{document}

\title{Hopf algebras of dimension $pq$ are trivial}
\author[X. Zhen]{Xiangjun Zhen$^\dag$}
\author[G. Liu]{Gongxiang Liu}
\author[J. Yu]{Jing Yu}
\address{Department of Mathematics, Nanjing University, Nanjing 210093, China}
\email{xjzhen@smail.nju.edu.cn}
\address{Department of Mathematics, Nanjing University, Nanjing 210093, China}
\email{gxliu@nju.edu.cn}
\address{School of Mathematical Sciences, University of Science and Technology of China, Hefei 230026, China}
\email{yujing46@ustc.edu.cn}
\thanks{$^\dag$ Corresponding author}
\date{}

\begin{abstract}
Let $H$ be a Hopf algebra of dimension $pq$ over an algebraically closed field $\k$ of characteristic zero, where $p$ and $q$ are distinct  primes.
We prove that $H$ is semisimple and hence, by the theorem of Etingof and Gelaki~\cite{EG98}, isomorphic to a group algebra or its dual.
\end{abstract}

\maketitle

\section{Introduction}

Classifying Hopf algebras of a given dimension is a basic problem in Hopf algebra theory (see the survey by Beattie and Garc\'{\i}a~\cite{BG13Survey}). Over an algebraically closed field $\k$ of characteristic zero, the classification is complete in dimensions $p$~\cite{Zhu94} and $p^2$~\cite{AS98,Mas96,Ng02} for every prime $p$, and in dimension $2p^2$~\cite{AN01,HN09,Nat99} for every odd prime $p$. In dimension $p^3$, the semisimple~\cite{Mas95Eight,Mas95Cube} and pointed~\cite{AS98Cube,CD98,SVO98} cases have been settled for every prime $p$, with a full classification for $p=2$~\cite{Stefan99} and $p=3$~\cite{BG13}. In dimension $pq^2$, with $p$ and $q$ distinct primes and $p>2$, see~\cite{ENO11,Nat99} for the semisimple case and~\cite{AN01} for the pointed case. Hopf algebras have also been classified in several low dimensions, such as $16$~\cite{CDR00,CDMM04,GV10,Kas00} and $30$~\cite{Fuk10,Nat07Memo}. It is worth noting that Ng and Wang~\cite{NW26} recently completed the classification of nonsemisimple Hopf algebras of dimension $24$.

We now turn to the classification of Hopf algebras of dimension $pq$, where $p$ and $q$ are distinct primes. Masuoka~\cite{Mas95Two} classified the semisimple Hopf algebras of dimension $2p$, and Gelaki and Westreich~\cite{GW00} obtained the classification in dimension $3p$ ($p>3$). Etingof and Gelaki~\cite{EG98} settled the general case and proved that every semisimple Hopf algebra of dimension $pq$ is isomorphic to $\k G$ or $(\k G)^*$ for some group $G$ of order $pq$. It has long been conjectured that every Hopf algebra of dimension $pq$ is semisimple, and despite a number of positive results for particular pairs of primes, the conjecture remained open in general (see~\cite{BG13Survey}). This is also a special case of Andruskiewitsch's  question on the semisimplicity of Hopf algebras of square-free dimension~\cite[Question~6.5]{And02}.
Once semisimplicity is established, the classification of Hopf algebras of dimension $pq$ is complete. 

Williams~\cite{Wil88} settled dimensions $6$ and $10$, and Andruskiewitsch and Natale~\cite{AN01} settled dimensions $15$, $21$, and $35$. Beattie and D\u{a}sc\u{a}lescu~\cite{BD04} proved semisimplicity in dimension $14$ and for  primes satisfying $2<p<q\leq\min\{3p+1,13\}$.
Fukuda~\cite{Fuk08} also established semisimplicity in several dimensions. Ng proved semisimplicity in dimension $2p$~\cite{Ng2p05} and in dimension $pq$ for twin primes $p$ and $q$~\cite{NgPQ04}. For odd primes $p<q$, Etingof and Gelaki~\cite{EG04} proved semisimplicity over $\mathbb{C}$ whenever $q\leq 2p+1$. Ng~\cite{NgPQII} extended this range to $p<q\leq 4p+11$.

In this paper, we prove that Hopf algebras of dimension $pq$ over $\k$ are semisimple when $p$ and $q$ are distinct odd primes. Together with Ng's theorem in dimension $2p$~\cite{Ng2p05}, this establishes the semisimplicity conjecture for all distinct primes $p$ and $q$ (see Theorem~\ref{thm:odd-pq-semisimple}). The classification of Hopf algebras of dimension $pq$ now follows from the theorem of Etingof and Gelaki~\cite{EG98}.

The organization of the paper is as follows. Section~\ref{sec:preliminaries} presents the definitions, results, and constructions used in the proof and concludes with an outline of the argument in Section~\ref{subsec:idea-of-proof}. In Section~\ref{sec:proof}, we prove the main theorem.

\section{Preliminaries}\label{sec:preliminaries}

Throughout this paper $\k$ denotes an algebraically closed field of characteristic 0 and all spaces are over $\k$. The set of nonnegative integers is denoted by $\mathbb{Z}_{+}$. The symbol $H$ will always denote a finite-dimensional Hopf algebra over $\k$.
Unless otherwise stated, all modules are understood to be left modules. Let $\operatorname{Rep}(H)$ denote the category of finite-dimensional left $H$-modules. 
For each $V\in\operatorname{Rep}(H)$, let $P(V)$ denote its projective cover.  Denote by $ G(H)$ the set of group-like elements of $ H $. 
We refer to \cite{MonS} for the basics about Hopf algebras.

\subsection{Grothendieck ring and character algebra}\label{subsec:grothendieck-ring}
Let $\operatorname{Gr}(H)$ be the Grothendieck ring of $\operatorname{Rep}(H)$ (see \cite[Definitions~1.5.8 and~4.5.2]{EGNO15}).
Its underlying additive group is the free abelian group generated by the isomorphism classes $[S_1],\ldots,[S_t]$ of the simple $H$-modules.
To every $V\in\operatorname{Rep}(H)$ we canonically associate its class
\[
   [V]=\sum_{i=1}^t [V:S_i]\,[S_i],
\]
where $[V:S_i]$ is the multiplicity of $S_i$ in a composition series of $V$.
For every short exact sequence $0\to U\to V\to W\to 0$ in $\operatorname{Rep}(H)$, we have $[V]=[U]+[W]$.
Its multiplication is defined by $[U][W]=[U\otimes W]$.

For $V\in\operatorname{Rep}(H)$ with action $\rho_V\colon H\to\operatorname{End}_{\k}(V)$, its character is the functional $\chi_V\in H^*$ given by $\chi_V(h)=\operatorname{Tr}_V(\rho_V(h))$ for $h\in H$. The \emph{character algebra} $R(H)$ (see~\cite[\S2.1]{L98}) is the $\k$-subalgebra of $H^*$ generated by these characters, with multiplication given by convolution.
In fact, we have the following isomorphism of $\k$-algebras (see \cite[Theorem~1.6(b) and Proposition~3.6]{L97})
\begin{equation}\label{eq:character-isomorphism}
   A:=\k\otimes_{\mathbb{Z}}\operatorname{Gr}(H)
   \xrightarrow{\ \sim\ } R(H),\qquad
   c\otimes[V]\longmapsto c\,\chi_V .
\end{equation}
In particular, characters are additive on short exact sequences and satisfy $\chi_{U\otimes V}=\chi_U * \chi_V$, where $*$ denotes convolution.
The irreducible characters $\chi_{S_1},\ldots,\chi_{S_t}$ form a $\k$-basis of $R(H)$.
When $H=\k G$ for a finite group $G$, this is a well-known result in character theory (see \cite[\S18.3, Theorem~15 and Proposition~17]{DF04}).
Via this isomorphism, we regard $A$ as a subalgebra of $H^*$.
Let $\varepsilon_{H^*}\colon H^*\to\k$ be the counit, given by $\varepsilon_{H^*}(f)=f(1_H)$. For every $V\in\operatorname{Rep}(H)$, we have
\begin{equation}\label{eq:counit-dimension}
    \varepsilon_{H^*}([V])=\chi_V(1_H)=\operatorname{Tr}_V(\operatorname{id}_V)=\dim_{\k}V.
\end{equation}
Thus $\varepsilon_{H^*}$ sends the class of each module to its dimension.

We conclude the subsection with the following lemma.

\begin{lemma}\label{lem:character-trace}
Let $H$ be an $n$-dimensional Hopf algebra over $\k$. For each $V\in\operatorname{Rep}(H)$, define $L_{\chi_V}\colon H^*\to H^*$ by $L_{\chi_V}(f)=\chi_V*f$. Then $\operatorname{Tr}_{H^*}(L_{\chi_V})\in n\mathbb{Z}$.
\end{lemma}

\begin{proof}
Consider the  character $\chi_{H^*}\in H^{**}$ of the regular $H^*$-module $H^*$. Let $i\colon H\to H^{**}$ be the canonical Hopf algebra isomorphism and let $z\in H$ be the unique element such that $i(z)=\frac{1}{n}\chi_{H^*}$. For every $f\in H^*$, we have
\[
f(z)=i(z)(f)=\frac{1}{n}\chi_{H^*}(f)=\frac{1}{n}\operatorname{Tr}_{H^*}(L_f),
\]
where $L_f$ denotes left multiplication by $f$ on $H^*$.
As observed in \cite[Section~3]{KMS22}, $z$ is idempotent. Taking $f=\chi_V$ gives
\[
\operatorname{Tr}_{H^*}(L_{\chi_V})=n\chi_V(z)=n\operatorname{Tr}_V(\rho_V(z))=n\dim_{\k}(zV). \qedhere
\]
\end{proof}

\subsection{Properties of Hopf algebras of dimension \texorpdfstring{$pq$}{pq}}\label{subsec:nonsemisimple-pq}
Since Ng~\cite[Theorem~3.3]{Ng2p05} proved that Hopf algebras of dimension $2q$ are semisimple for every odd prime $q$, we assume throughout the rest of the paper that $p<q$ are odd primes unless otherwise stated.
We begin with two results of Ng~\cite{NgPQII}.
For $\beta\in G(H^*)$, let $\k_\beta$ denote the one-dimensional $H$-module on which each $h\in H$ acts by the scalar $\beta(h)$.
The following lemma is an analogue of \cite[Lemma~1.4]{NgPQII} and follows by taking left duals.

\begin{lemma}[{\cite[Lemma~1.4]{NgPQII}}]\label{lem:right-tensor-divisibility}
Let $H$ be a finite-dimensional Hopf algebra over $\k$, let $\beta\in G(H^*)$, and let $V\in\operatorname{Rep}(H)$ be indecomposable.
If $V\otimes\k_\beta\cong V$ as $H$-modules, then $\operatorname{ord}(\beta)\mid\dim V$.
\end{lemma}

The following result  is crucial for deriving the contradiction in the proof of our main theorem.

\begin{lemma}[{\cite[Lemma~2.1(vi)]{NgPQII}}]\label{lem:ng-projective-cover}
Let $H$ be a nonsemisimple Hopf algebra of dimension $pq$ over $\k$, where $p<q$ are odd primes.
Then there exists a simple $H$-module $V$ such that
\[
    \dim V \geq 2 \qquad\text{and}\qquad \dim P(V) \geq 2\dim V.
\]
\end{lemma}

We now turn to the Hopf algebras of dimension $pq$ studied in this paper.
Suppose that $H$ is a nonsemisimple Hopf algebra of dimension $pq$ over $\k$.
By \cite{LR88}, $H^*$ is also nonsemisimple.
By Ng~\cite[Lemma~2.1(i)]{NgPQII}, at least one of $H$ and $H^*$ is not unimodular.
Thus, after replacing $H$ by $H^*$ if necessary, we may assume that $H$ is not unimodular.
The distinguished group-like element $\alpha\in G(H^*)$ is therefore different from the counit $\varepsilon$ of $H$, so $G(H^*)$ is nontrivial.
Applying \cite[Lemma~2.1(ii)]{NgPQII} to $H^*$, we obtain $|G(H^*)|=p$ and $G(H^*)=\langle\alpha\rangle$.

In the rest of the paper, except in Section~\ref{subsec:idea-of-proof}, we keep $H$ fixed as above. We regard $A=\k\otimes_{\mathbb{Z}}\operatorname{Gr}(H)$ as a subalgebra of $H^*$ via the isomorphism~\eqref{eq:character-isomorphism}.

\subsection{Actions on \texorpdfstring{$\operatorname{Irr}(H)$}{Irr(H)}}\label{subsec:actions-on-irr}
Let $\operatorname{Irr}(H)$ denote the set of isomorphism classes of simple $H$-modules, and write $[V]$ for the class of a simple $H$-module $V$.
This notation agrees with that of Section~\ref{subsec:grothendieck-ring}, where $[V]$ denotes the corresponding basis element of $\operatorname{Gr}(H)$.

Following Ng~\cite[Section~3]{NgPQII}, we consider the left and right actions of $G(H^*)$ on $\operatorname{Irr}(H)$.
For $[V]\in\operatorname{Irr}(H)$ and $\beta\in G(H^*)$, these actions are given by
\[
    \beta\cdot[V]=[\k_\beta\otimes V],
    \qquad
    [V]\cdot\beta=[V\otimes\k_\beta].
\]
The \emph{left orbit} $O_l(V)$ and the \emph{right orbit} $O_r(V)$ of $V$ are the orbits of $[V]$ under the left and right actions, respectively.
In Ng's terminology, $V$ is \emph{regular} if $O_l(V)=O_r(V)$.
Since the tensor product is associative and $G(H^*)$ is abelian, we also have a left action of $G(H^*)\times G(H^*)$ on $\operatorname{Irr}(H)$.
For $[V]\in\operatorname{Irr}(H)$ and $\beta,\gamma\in G(H^*)$, this action is given by
\[
    (\beta,\gamma)\cdot[V]=[\k_\beta\otimes V\otimes\k_\gamma].
\]
Write $O(V)$ for the orbit of $[V]$ under this action.
Part~(ii) of the following lemma follows from~\cite[Remark~3.1]{NgPQII}.

\begin{lemma}\label{lem:joint-orbit-properties}
Let $V$ be a simple $H$-module.
\begin{enumerate}[(i)]
    \item We have
    \begin{equation}\label{eq:joint-orbit-decomposition}
        O(V)=\bigcup_{\beta\in G(H^*)} O_l(V\otimes\k_\beta)
            =\bigcup_{\beta\in G(H^*)} O_r(\k_\beta\otimes V).
    \end{equation}
    \item If $[W]\in O(V)$, then
    \begin{equation}\label{eq:orbit-projective-dimension}
        \dim_{\k}W=\dim_{\k}V,\qquad \dim_{\k}P(W)=\dim_{\k}P(V).
    \end{equation}
\end{enumerate}
\end{lemma}

Under the additional assumption $q\leq 4p+11$, Ng shows that every simple $H$-module $V$ is regular and $|O_l(V)|=p$ (see \cite[Lemmas~4.2--4.4]{NgPQII}). In our notation,
\[
    O(V)=O_l(V)=O_r(V)
    =\bigl\{[V],[\k_\alpha\otimes V],\ldots,[\k_{\alpha^{p-1}}\otimes V]\bigr\}.
\]
Ng then groups the terms in the dimension formula for $H$ by $G(H^*)$-orbit (see the formula following Lemma~4.8 in \cite{NgPQII}) and derives $q>4p+11$, a contradiction.
For arbitrary odd primes $p<q$, we pursue a similar counting argument (see~\eqref{eq:orbit-dimension-formula}) by decomposing each $O(V)$ into its distinct right orbits.

\subsection{Orbit decompositions and the associated multiplication operator}\label{subsec:orbit-decompositions}

We now assume that the action of $G(H^*)\times G(H^*)$ on $\operatorname{Irr}(H)$ has $r+1$ orbits $O(X_0),O(X_1),\ldots,O(X_r)$.
We take $X_0=\k_\alpha$. 
For each $0\leq i\leq r$, let $m_i$ be the number of distinct right orbits in the decomposition~\eqref{eq:joint-orbit-decomposition} of $O(X_i)$ and choose one representative $[X_{i,j}]$ from each, so that
\[
    \operatorname{Irr}(H)=\bigsqcup_{i=0}^{r} O(X_i)
    =\bigsqcup_{i=0}^{r}\bigsqcup_{j=1}^{m_i} O_r(X_{i,j}).
\]

Fix $0\leq i\leq r$ and $1\leq j\leq m_i$. By~\eqref{eq:joint-orbit-decomposition}, we may assume that $X_{i,j}=\k_\beta\otimes X_i$ for some $\beta\in G(H^*)$.
Since $G(H^*)=\langle\alpha\rangle$ has prime order $p$, every right orbit has size $1$ or $p$. We have
\[
    |O_r(X_{i,j})|=1 \iff (\k_\beta\otimes X_i)\otimes\k_\alpha\cong\k_\beta\otimes X_i \iff X_i\otimes\k_\alpha\cong X_i,
\]
where the last equivalence follows by tensoring on the left with $\k_{\beta^{-1}}$.
Thus all right orbits in $O(X_i)$ have the same size, which we denote by $s_i$. In particular,
\begin{equation}\label{eq:right-orbit-size}
    s_i=|O_r(X_i)|\in\{1,p\},\qquad
    s_i=1 \iff X_i\otimes\k_\alpha\cong X_i.
\end{equation}
It follows that $|O(X_i)|=m_i s_i$ for $0\leq i\leq r$.
For $X_0=\k_\alpha$, we have $m_0=1$ and $s_0=p$, so $|O(X_0)|=p$.
For $0\leq i\leq r$, set
\[
    d_i=\dim_{\k}X_i,\qquad h_i=\dim_{\k}P(X_i).
\]
Then $h_i\geq d_i\geq 2$ for $1\leq i\leq r$.
The dimension formula for the regular $H$-module (see (3.11) in \cite{NgPQII}), together with~\eqref{eq:orbit-projective-dimension}, gives the following lemma.

\begin{lemma}\label{lem:orbit-dimension-formula}
We have
\begin{equation}\label{eq:orbit-dimension-formula}
    pq=p\dim_{\k}P(\k)+\sum_{i=1}^{r}m_i s_i d_i h_i.
\end{equation}
\end{lemma}

The dimension formula~\eqref{eq:orbit-dimension-formula} is central to our proof.
We now introduce some notation and record a few results for the proof.
We use uppercase letters for objects of $\operatorname{Rep}(H)$ and lowercase letters for elements of the character algebra $A$.
For each orbit $O(X_i)$, set
\[
    V_i=\bigoplus_{[V]\in O(X_i)} V\in\operatorname{Rep}(H),
    \qquad
    v_i=[V_i]=\sum_{[V]\in O(X_i)}[V]\in A.
\]
The isomorphism~\eqref{eq:character-isomorphism} sends $[\k_\alpha]$ to $\alpha$, so $\alpha$ is a character. Thus $G(H^*)=\langle\alpha\rangle\subseteq A$. Set
\[
    e=\frac{1}{p}\sum_{i=0}^{p-1}\alpha^i.
\]
Then $e$ is an idempotent in $A$ and satisfies $e\alpha=\alpha e=e$. Also note that $pe=v_0$.

Let $[X]\in\operatorname{Irr}(H)$ and assume that $[X]\in O(X_i)$ for some $0\leq i\leq r$. Then
\[
    \begin{aligned}
        e[X]e
        &=\frac{1}{p^2}\sum_{j,k=0}^{p-1}[\k_{\alpha^j}\otimes X\otimes\k_{\alpha^k}]\\
        &=\frac{1}{p^2}\sum_{(\beta,\gamma)\in G(H^*)\times G(H^*)}(\beta,\gamma)\cdot[X]\\
        &=\frac{1}{|O(X_i)|}\sum_{[V]\in O(X_i)}[V]\\
        &=\frac{1}{m_i s_i}v_i.
    \end{aligned}
\]
The penultimate equality holds because each class in $O(X_i)$ occurs exactly $p^2/|O(X_i)|$ times in the sum over $G(H^*)\times G(H^*)$ (see~\cite[\S4.1, Proposition~2]{DF04}).
Set
\[
    B=eAe\subseteq A.
\]
Since the irreducible characters form a basis of $A$ (see Section~\ref{subsec:grothendieck-ring}), the preceding computation shows that $\{v_0,\ldots,v_r\}$ is a $\k$-basis of $B$, so $\dim_{\k} B=r+1$.
Also note that $B=\{a\in A\mid ea=a=ae\}$.
Set
\[
    E=eH^*.
\]
Then $B\subseteq E$, since $eb=b$ for all $b\in B$.
By the Nichols--Zoeller theorem~\cite{NZ89}, $H^*$ is a free $\k G(H^*)$-module of rank $q$.
Note that $e\,\k G(H^*)=\k e$, so $\dim_{\k} E=q$.

Finally, for each $0\leq i\leq r$, we use the chosen representatives of the right orbits in $O(X_i)$ to define
\[
    Y_i=\bigoplus_{j=1}^{m_i}X_{i,j}\in\operatorname{Rep}(H),
    \qquad y_i=[Y_i]\in A.
\]
For all $0\leq i\leq r$, we have $y_i e=\frac{1}{s_i}v_i\in B$. Also note that $\dim_{\k}Y_i=m_i d_i$.
Choose $1\leq k\leq r$ such that $\dim_{\k}Y_k=\min_{1\leq i\leq r}\dim_{\k}Y_i$, and set
\begin{equation}\label{eq:definition-of-b}
    Y=Y_k,\qquad N=\dim_{\k}Y,\qquad y=[Y],\qquad b=ye=\frac{1}{s_k}v_k\in B.
\end{equation}
Let $L_b\colon H^*\to H^*$ denote left multiplication by $b$, given by $L_b(f)=bf$.
Since $b\in eAe$, $L_b$ preserves both $B$ and $E$. The characteristic polynomials of $L_b|_B$ and $L_b|_E$ will play a central role in the proof.

\subsection{Idea of the proof}\label{subsec:idea-of-proof}
Here we outline the main idea of our proof. This will help clarify the constructions in the preceding subsection.

The constructions in Sections~\ref{subsec:actions-on-irr} and~\ref{subsec:orbit-decompositions} apply to any Hopf algebra $H$ of dimension $pq$ with $|G(H^*)|=p$, regardless of whether $H$ is semisimple.
\textbf{\boldmath Assume for the moment that $H$ is a semisimple Hopf algebra of dimension $pq$ with $|G(H^*)|=p$.}
In this subsection, let $\alpha$ be a generator of $G(H^*)$.
We now describe these constructions explicitly in the semisimple case.

By the theorem of Etingof and Gelaki~\cite{EG98}, $H\cong\k G$ or $H\cong(\k G)^*$ for some group $G$ of order $pq$.
We recall the classification of groups of order $pq$ from~\cite[\S5.5]{DF04}. If $p\nmid(q-1)$, then $G$ is cyclic. In this case $(\k G)^*\cong\k G$ as Hopf algebras, so $|G(H^*)|=pq$, contrary to our assumption.
If $p\mid(q-1)$, then $G$ is either cyclic or isomorphic to the unique nonabelian semidirect product $C_q\rtimes C_p$. The cyclic case is excluded as above. Moreover, if $H\cong(\k G)^*$, then $H^*\cong\k G$ has $pq$ group-like elements.
Thus $H\cong\k G$ with $G=C_q\rtimes C_p$. Indeed, $G/[G,G]\cong C_p$, so $G(H^*)\cong\operatorname{Hom}(G,\k^\times)$ has order $p$.

Up to isomorphism, $H$ has $p$ one-dimensional simple modules and $(q-1)/p$ $p$-dimensional simple modules (see~\cite[\S19.3, Proposition~13(2)--(3)]{DF04}).
Let $\beta\in G(H^*)$, and let $V$ be any $p$-dimensional simple $H$-module. As shown in the proof of that proposition, $\chi_{\k_\beta}$ equals $1$ on $C_q$ and $\chi_V$ vanishes on $G\setminus C_q$.
It follows that $\chi_{\k_\beta}*\chi_V=\chi_V=\chi_V*\chi_{\k_\beta}$, and hence
\[
    \k_\beta\otimes V\cong V\cong V\otimes\k_\beta.
\]

We now apply the constructions of Section~\ref{subsec:orbit-decompositions} to $H$. By the preceding isomorphisms, $r=(q-1)/p$ and, for every $1\leq i\leq r$,
\[
    O(X_i)=O_l(X_i)=O_r(X_i)=\{[X_i]\}.
\]
It follows that $m_i=s_i=1$ for $1\leq i\leq r$. Since $H$ is semisimple, $P(X_i)\cong X_i$ for all $0\leq i\leq r$. Hence $h_0=d_0=1$ and $d_i=h_i=p$ for $1\leq i\leq r$.
Then~\eqref{eq:orbit-dimension-formula} becomes $pq=p+rp^2$. This identity also follows directly from $r=(q-1)/p$.
Dividing by $p$, we obtain
\[
    q=1+rp\leq (1+r)p.
\]
For $1\leq i\leq r$, we have $V_i=Y_i=X_i$. Thus $N=p$, and any of $X_1,\ldots,X_r$ may be taken as $Y$. The preceding inequality therefore reads
\[
    q\leq (1+r)N.
\]

Motivated by the semisimple case, we will prove the same bound in the next section for the Hopf algebra $H$ considered in Sections~\ref{subsec:actions-on-irr} and~\ref{subsec:orbit-decompositions} (see Lemma~\ref{lem:trace-congruence}). On the other hand, Lemma~\ref{lem:ng-projective-cover} and the dimension formula~\eqref{eq:orbit-dimension-formula} give $q>(1+r)N$ (see Lemma~\ref{lem:strict-lower-bound}). These bounds are incompatible, so no such $H$ can exist. The main theorem then follows from the reduction in Section~\ref{subsec:nonsemisimple-pq}.

\section{Proof of the main theorem}\label{sec:proof}
We retain the notation of Sections~\ref{subsec:nonsemisimple-pq}--\ref{subsec:orbit-decompositions}.
Recall that $p<q$ are odd primes and that $H$ is a nonsemisimple Hopf algebra of dimension $pq$ with $|G(H^*)|=p$.
We will prove both $q>(1+r)N$ and $q\leq (1+r)N$. This contradiction shows that no such $H$ exists.

\subsection{\texorpdfstring{{\boldmath$q>(1+r)N$}}{q > (1+r)N}}
The lower bound follows readily from Lemma~\ref{lem:ng-projective-cover} and the dimension formula~\eqref{eq:orbit-dimension-formula}.

\begin{lemma}\label{lem:strict-lower-bound}
We have $q>(1+r)N$.
\end{lemma}

\begin{proof}
By~\eqref{eq:orbit-dimension-formula}, we have
\[
    q=\dim_{\k}P(\k)+\sum_{i=1}^{r}m_i d_i\frac{s_i h_i}{p}.
\]
Recall that $N=\min_{1\leq i\leq r}m_i d_i$. We first show that for every $1\leq i\leq r$,
\[
    m_i d_i\frac{s_i h_i}{p}\geq N.
\]
If $s_i=p$, then $h_i\geq d_i\geq 2$, so
\[
    m_i d_i\frac{s_i h_i}{p}=m_i d_i h_i\geq 2m_i d_i\geq 2N.
\]
If $s_i=1$, then $X_i\otimes\k_\alpha\cong X_i$ by~\eqref{eq:right-orbit-size}. Lemma~\ref{lem:right-tensor-divisibility} gives $p\mid d_i$. Hence
\[
    m_i d_i\frac{s_i h_i}{p}=m_i h_i\frac{d_i}{p}\geq m_i h_i\geq m_i d_i\geq N.
\]
Moreover, if $h_i\geq 2d_i$, then $m_i h_i\geq 2m_i d_i\geq 2N$ in this case as well.

By Lemma~\ref{lem:ng-projective-cover}, there is an index $1\leq i\leq r$ such that $h_i\geq 2d_i$. For this $i$, the preceding estimates give $m_i d_i\frac{s_i h_i}{p}\geq 2N$. Therefore
\[
    q\geq \dim_{\k}P(\k)+(r-1)N+2N
      =\dim_{\k}P(\k)+(1+r)N>(1+r)N.\qedhere
\]
\end{proof}

\subsection{\texorpdfstring{{\boldmath$q\leq (1+r)N$}}{q ≤ (1+r)N}}
The proof of the upper bound is more involved. We first estimate the trace of $L_b|_B$.  Set
\[
    M=\operatorname{Tr}_B\bigl(L_b|_B\bigr).
\]
We  show that $M$ is a nonnegative integer and that $(1+r)N-M>0$.

\begin{lemma}\label{lem:trace-bounds}
We have $M\in\mathbb{Z}$ and $0\leq M\leq rN$.
\end{lemma}

\begin{proof}
We compute the trace of $L_b|_B$ in the basis $\{v_0,\ldots,v_r\}$. By~\eqref{eq:definition-of-b}, $b=ye=\frac{1}{s_k}v_k$. Since $v_0=pe$, we have
\[
    L_b(v_0)=bv_0=pbe=pb=\frac{p}{s_k}v_k.
\]
As $k\geq 1$, the coefficient of $v_0$ is zero.

For $1\leq j\leq r$, write
\begin{equation}\label{eq:bvj-expansion}
    bv_j=\sum_{i=0}^{r}M_{ij}v_i=\sum_{i=0}^{r}M_{ij}[V_i],
\end{equation}
where $M_{ij}\in\k$.
Since $v_j\in B$, we have $ev_j=v_j$, and hence $bv_j=yev_j=yv_j=[Y\otimes V_j]$. Thus both sides of~\eqref{eq:bvj-expansion} lie in $\operatorname{Gr}(H)$. Comparing the coefficients of $[X_i]$ gives
\[
    M_{ij}=[Y\otimes V_j:X_i]\in\mathbb{Z}_{+}.
\]
Thus the trace is a nonnegative integer given by
\[
    M=\sum_{j=1}^{r}M_{jj}=\sum_{j=1}^{r}[Y\otimes V_j:X_j].
\]
Applying $\varepsilon_{H^*}$ to~\eqref{eq:bvj-expansion} and using~\eqref{eq:counit-dimension}, we obtain
\[
    N\dim_{\k}V_j=\dim_{\k}(Y\otimes V_j)
    =\sum_{i=0}^{r}M_{ij}\dim_{\k}V_i\geq M_{jj}\dim_{\k}V_j.
\]
Since $\dim_{\k}V_j>0$, we have $M_{jj}\leq N$ for each $1\leq j\leq r$. Consequently, $0\leq M\leq rN$.
\end{proof}

We next prove that $(1+r)N\equiv M\pmod q$. This congruence implies that $q$ divides $(1+r)N-M$. By Lemma~\ref{lem:trace-bounds}, this difference is a positive integer, so $q\leq (1+r)N-M$. To prove the congruence, set
\[
    T=L_b|_E\colon E\to E.
\]
We begin with the following result.

\begin{lemma}\label{lem:power-traces}
For every integer $i\geq 1$, we have $\operatorname{Tr}_E(T^i)\in q\mathbb{Z}$.
\end{lemma}

\begin{proof}
Fix $i\geq 1$ and set $W=Y^{\otimes i}\otimes V_0$. Since $b=ye$ and $eb=b$, we have $b^i=y^i e$. Under the isomorphism~\eqref{eq:character-isomorphism}, we  have
\[
    \chi_W=y^i v_0=p y^i e=p b^i.
\]
Thus $L_{\chi_W}=pL_{b^i}$. Since $e$ is idempotent, $H^*=E\oplus(1-e)H^*$, where $E=eH^*$. On $E$, left multiplication by $b^i$ is $(L_b|_E)^i=T^i$. On $(1-e)H^*$ it is zero, since $b^i(1-e)=y^i e(1-e)=0$. Taking traces gives
\[
    \operatorname{Tr}_{H^*}(L_{\chi_W})
    =p\operatorname{Tr}_{H^*}(L_{b^i})
    =p\operatorname{Tr}_E(T^i).
\]
By Lemma~\ref{lem:character-trace}, the left-hand side is a multiple of $pq$. Hence $\operatorname{Tr}_E(T^i)\in q\mathbb{Z}$.
\end{proof}

Since $\k$ has characteristic zero, we regard $\mathbb{Q}$ as a subfield of $\k$. Let $\mathbb{Z}_{(q)}$ denote the subring of $\mathbb{Q}$ consisting of fractions $a/b$ in lowest terms such that $\gcd(b,q)=1$.
Let $\mathbb{F}_q=\mathbb{Z}/q\mathbb{Z}$. For each $a\in\mathbb{Z}$, write $\overline{a}$ for its image in $\mathbb{F}_q$ under the quotient map.
For $x\in\mathbb{Z}_{(q)}$, choose any representation $x=a/b$ with $a,b\in\mathbb{Z}$ and $\gcd(b,q)=1$, and set
\[
    \overline{x}=\overline{a}\,\overline{b}^{-1}\in\mathbb{F}_q.
\]
The map $\phi\colon\mathbb{Z}_{(q)}\to\mathbb{F}_q$, $x\mapsto\overline{x}$, is a well-defined surjective ring homomorphism. Clearly, $\phi(q\mathbb{Z}_{(q)})=0$.
This induces a ring homomorphism
\[
    \phi\colon \mathbb{Z}_{(q)}[t]\longrightarrow\mathbb{F}_q[t],\qquad
    \sum_{i=0}^{n}a_i t^i\longmapsto\sum_{i=0}^{n}\phi(a_i)t^i.
\]
Let $F(t)\in\k[t]$ denote the characteristic polynomial of $T$. We have the following result.

\begin{lemma}\label{lem:reduced-charpoly}
We have $F(t)\in\mathbb{Z}_{(q)}[t]$ and $\phi(F(t))=(t-\overline{N})^q$.
\end{lemma}

\begin{proof}
Since $\dim_{\k}E=q$, write
\[
    F(t)=t^q+c_1t^{q-1}+\cdots+c_q,
\]
and set $c_0=1$. For $1\leq j<q$, Newton's formulas (see~\cite[\S14.6, Exercise~22]{DF04}) give
\[
    j c_j=-\sum_{i=1}^{j}c_{j-i}\operatorname{Tr}_E(T^i).
\]
We show by induction that $c_j\in q\mathbb{Z}_{(q)}$ for $1\leq j<q$. Since $c_0=1$, we may assume that $c_0,\ldots,c_{j-1}\in\mathbb{Z}_{(q)}$. By Lemma~\ref{lem:power-traces}, each summand on the right lies in $q\mathbb{Z}_{(q)}$. As $1\leq j<q$ and $q$ is prime, $j$ is invertible in $\mathbb{Z}_{(q)}$. Hence $c_j\in q\mathbb{Z}_{(q)}$, as claimed.

To determine $c_q$, take a nonzero left integral $\lambda\in H^*$, which exists since $H^*$ is finite dimensional (see~\cite[Chapter~2]{MonS}). Then $f\lambda=\varepsilon_{H^*}(f)\lambda$ for every $f\in H^*$. Since $pe=v_0=[V_0]$ and $\dim_{\k}V_0=p$, equation~\eqref{eq:counit-dimension} gives $\varepsilon_{H^*}(e)=1$. Also, $b=ye$ gives $\varepsilon_{H^*}(b)=\dim_{\k}Y=N$. Thus
\[
    e\lambda=\lambda,\qquad b\lambda=N\lambda.
\]
The first equality gives $\lambda\in eH^*=E$, and the second gives $T(\lambda)=N\lambda$. Hence $N$ is an eigenvalue of $T$ and $F(N)=0$. It follows that
\[
    c_q=-N^q-\sum_{j=1}^{q-1}c_j N^{q-j}\in\mathbb{Z}_{(q)}.
\]
Thus $F(t)\in\mathbb{Z}_{(q)}[t]$. Since $\phi(c_j)=0$ for $1\leq j<q$, the formula above gives $\phi(c_q)=-\overline{N}^{\,q}$. Hence we have
\[
    \phi(F(t))=t^q-\overline{N}^{\,q}=(t-\overline{N})^q.\qedhere
\]
\end{proof}

Let $G(t)\in\k[t]$ denote the characteristic polynomial of $L_b|_B$. The coefficient of $t^r$ in $G(t)$ will give the desired congruence $(1+r)N\equiv M\pmod q$.

\begin{lemma}\label{lem:trace-congruence}
We have $(1+r)N\equiv M\pmod q$. In particular, $q\leq (1+r)N$.
\end{lemma}

\begin{proof}
It suffices to prove that $(1+r)N\equiv M\pmod q$, since the inequality then follows as explained before Lemma~\ref{lem:power-traces}.

Write the matrix of $L_b|_B$ with respect to the basis $(v_0,\ldots,v_r)$ of $B$ as $M'=(M_{ij})_{0\leq i,j\leq r}$. By the proof of Lemma~\ref{lem:trace-bounds}, its entries are integers. Thus $G(t)$ is a monic polynomial in $\mathbb{Z}[t]$ of degree $r+1$.
Extend $(v_0,\ldots,v_r)$ to an ordered basis $\mathcal{B}=(v_0,\ldots,v_r,w_1,\ldots,w_s)$ of $E$. With respect to this basis, the matrix of $T$ has the form
\[
    [T]_{\mathcal{B}}=
    \begin{pmatrix}
        M' & C\\
        0 & D
    \end{pmatrix}.
\]
Then $F(t)=G(t)K(t)$, where $K(t)=\det(tI_s-D)$.

By Lemma~\ref{lem:reduced-charpoly}, $F(t)\in\mathbb{Z}_{(q)}[t]$. Since $G(t)\in\mathbb{Z}[t]\subseteq\mathbb{Z}_{(q)}[t]$ is monic, it follows from~\cite[Proposition~11.2.9, p.~327]{Artin11} that $K(t)\in\mathbb{Z}_{(q)}[t]$. Applying $\phi$ to $F(t)=G(t)K(t)$ shows that $\phi(G(t))$ divides $\phi(F(t))=(t-\overline{N})^q$ in $\mathbb{F}_q[t]$. As $\phi(G(t))$ is monic of degree $r+1$, it follows that
\begin{equation}\label{eq:reduced-G-power}
    \phi(G(t))=(t-\overline{N})^{r+1}.
\end{equation}
On the other hand, the coefficient of $t^r$ in $G(t)$ is $-\operatorname{Tr}_B(L_b|_B)=-M$. Hence
\begin{equation}\label{eq:reduced-G-expansion}
    \phi(G(t))=t^{r+1}-\overline{M}t^r+\cdots.
\end{equation}
Comparing the coefficients of $t^r$ in~\eqref{eq:reduced-G-power} and~\eqref{eq:reduced-G-expansion} gives $\overline{M}=(1+r)\overline{N}$ in $\mathbb{F}_q$. Thus $(1+r)N\equiv M\pmod q$.
\end{proof}

\subsection{The main theorem}\label{subsec:main-theorem}

We can now prove our main theorem.

\begin{theorem}\label{thm:odd-pq-semisimple}
Let $p<q$ be  primes, and let $H$ be a Hopf algebra of dimension $pq$ over an algebraically closed field $\k$ of characteristic zero. Then $H$ is semisimple.
\end{theorem}

\begin{proof}
If $p=2$, then $H$ is semisimple by Ng's theorem in dimension $2q$~\cite[Theorem~3.3]{Ng2p05}. We may therefore assume that $p$ and $q$ are odd.
Suppose that $H$ is nonsemisimple. By the reduction in Section~\ref{subsec:nonsemisimple-pq}, after replacing $H$ by $H^*$ if necessary, we may assume that $G(H^*)=\langle\alpha\rangle$ has order $p$. Lemma~\ref{lem:strict-lower-bound} gives $q>(1+r)N$, while Lemma~\ref{lem:trace-congruence} gives $q\leq (1+r)N$. This contradiction shows that $H$ must be semisimple.
\end{proof}

By Etingof and Gelaki~\cite{EG98}, every such Hopf algebra is therefore isomorphic to $\k G$ or $(\k G)^*$ for some group $G$ of order $pq$. Recall from Section~\ref{subsec:idea-of-proof} that $C_{pq}$ is the only group of order $pq$ up to isomorphism if $p\nmid(q-1)$, while if $p\mid(q-1)$ there is also a unique nonabelian group $C_q\rtimes C_p$ up to isomorphism. Since $(\k C_{pq})^*\cong\k C_{pq}$ as Hopf algebras, we obtain the following classification.

\begin{theorem}\label{thm:pq-classification}
Let $p<q$ be primes and let $\k$ be an algebraically closed field of characteristic zero. The Hopf algebras of dimension $pq$ over $\k$ are as follows.
\begin{enumerate}
    \item If $p\nmid(q-1)$, the only such Hopf algebra is $\k C_{pq}$.
    \item If $p\mid(q-1)$, there are exactly three such Hopf algebras, namely
    \[
        \k C_{pq},\qquad \k G,\qquad (\k G)^*,
    \]
    where $G$ is the nonabelian semidirect product $C_q\rtimes C_p$.
\end{enumerate}
\end{theorem}

\section*{Acknowledgements}
The second author was supported by  National Key R\&D Program of China 2024YFA1013802 and NSFC 12271243.  The third author was supported by the China Postdoctoral Science Foundation 2026M793410 and the Postdoctoral Fellowship Program of China Postdoctoral Science Foundation GZC20261673.
The argument establishing $q\leq (1+r)N$ is due to ChatGPT 6 Astra. ChatGPT also helped refine the language of the manuscript. The authors simplified the argument and then edited and checked the manuscript. The authors take full responsibility for the content of this paper.


\begin{thebibliography}{99}
\hypersetup{urlcolor=black}
\setlength{\itemsep}{0em}
\bibitem{And02}N. Andruskiewitsch, About finite dimensional Hopf algebras, in: Quantum symmetries in theoretical physics and mathematics (Bariloche, 2000), Contemp. Math., vol. 294, Amer. Math. Soc., Providence, RI, 2002, pp. 1--57. \url{https://doi.org/10.1090/conm/294/04969}.
\bibitem{AN01}N. Andruskiewitsch, S. Natale, Counting arguments for Hopf algebras of low dimension, Tsukuba J. Math. 25 (1) (2001) 187--201. \url{https://tsukuba.repo.nii.ac.jp/records/17269}.
\bibitem{AS98}N. Andruskiewitsch, H.-J. Schneider, Hopf algebras of order $p^2$ and braided Hopf algebras of order $p$, J. Algebra 199 (2) (1998) 430--454. \url{https://doi.org/10.1006/jabr.1997.7175}.
\bibitem{AS98Cube}N. Andruskiewitsch, H.-J. Schneider, Lifting of quantum linear spaces and pointed Hopf algebras of order $p^3$, J. Algebra 209 (2) (1998) 658--691. \url{https://doi.org/10.1006/jabr.1998.7643}.
\bibitem{Artin11}M. Artin, Algebra, 2nd ed., Prentice Hall, Boston, MA, 2011.
\bibitem{BD04}M. Beattie, S. D\u{a}sc\u{a}lescu, Hopf algebras of dimension $14$, J. London Math. Soc. (2) 69 (1) (2004) 65--78. \url{https://doi.org/10.1112/S0024610703004927}.
\bibitem{BG13}M. Beattie, G. A. Garc\'{\i}a, Techniques for classifying Hopf algebras and applications to dimension $p^3$, Comm. Algebra 41 (8) (2013) 3108--3129. \url{https://doi.org/10.1080/00927872.2012.692004}.
\bibitem{BG13Survey}M. Beattie, G. A. Garc\'{\i}a, Classifying Hopf algebras of a given dimension, in: Hopf algebras and tensor categories, Contemp. Math., vol. 585, Amer. Math. Soc., Providence, RI, 2013, pp. 125--152. \url{https://doi.org/10.1090/conm/585/11615}.
\bibitem{CD98}S. Caenepeel, S. D\u{a}sc\u{a}lescu, Pointed Hopf algebras of dimension $p^3$, J. Algebra 209 (2) (1998) 622--634. \url{https://doi.org/10.1006/jabr.1998.7504}.
\bibitem{CDR00}S. Caenepeel, S. D\u{a}sc\u{a}lescu, \c{S}. Raianu, Classifying pointed Hopf algebras of dimension 16, Comm. Algebra 28 (2) (2000) 541--568. \url{https://doi.org/10.1080/00927870008826843}.
\bibitem{CDMM04}C. C\u{a}linescu, S. D\u{a}sc\u{a}lescu, A. Masuoka, C. Menini, Quantum lines over non-cocommutative cosemisimple Hopf algebras, J. Algebra 273 (2) (2004) 753--779. \url{https://doi.org/10.1016/j.jalgebra.2003.08.006}.
\bibitem{DF04}D. S. Dummit, R. M. Foote, Abstract algebra, 3rd ed., John Wiley \& Sons, Hoboken, NJ, 2004.
\bibitem{EG98}P. Etingof, S. Gelaki, Semisimple Hopf algebras of dimension $pq$ are trivial, J. Algebra 210 (2) (1998) 664--669. \url{https://doi.org/10.1006/jabr.1998.7568}.
\bibitem{EG04}P. Etingof, S. Gelaki, On Hopf algebras of dimension $pq$, J. Algebra 277 (2) (2004) 668--674. \url{https://doi.org/10.1016/j.jalgebra.2003.06.003}.
\bibitem{EGNO15}P. Etingof, S. Gelaki, D. Nikshych, V. Ostrik, Tensor categories, Mathematical Surveys and Monographs, vol. 205, American Mathematical Society, Providence, RI, 2015. \url{https://doi.org/10.1090/surv/205}.
\bibitem{ENO11}P. Etingof, D. Nikshych, V. Ostrik, Weakly group-theoretical and solvable fusion categories, Adv. Math. 226 (1) (2011) 176--205. \url{https://doi.org/10.1016/j.aim.2010.06.009}.
\bibitem{Fuk08}D. Fukuda, Structure of coradical filtration and its application to Hopf algebras of dimension $pq$, Glasg. Math. J. 50 (2) (2008) 183--190. \url{https://doi.org/10.1017/S0017089508004126}.
\bibitem{Fuk10}D. Fukuda, Hopf algebras of dimension 30, J. Algebra Appl. 9 (1) (2010) 11--15. \url{https://doi.org/10.1142/S0219498810003902}.
\bibitem{GV10}G. A. Garc\'{\i}a, C. Vay, Hopf algebras of dimension 16, Algebr. Represent. Theory 13 (4) (2010) 383--405. \url{https://doi.org/10.1007/s10468-009-9128-7}.
\bibitem{GW00}S. Gelaki, S. Westreich, On semisimple Hopf algebras of dimension $pq$, Proc. Amer. Math. Soc. 128 (1) (2000) 39--47. \url{https://doi.org/10.1090/S0002-9939-99-04961-8}.
\bibitem{HN09}M. Hilgemann, S.-H. Ng, Hopf algebras of dimension $2p^2$, J. London Math. Soc. (2) 80 (2) (2009) 295--310. \url{https://doi.org/10.1112/jlms/jdp026}.
\bibitem{Kas00}Y. Kashina, Classification of semisimple Hopf algebras of dimension 16, J. Algebra 232 (2) (2000) 617--663. \url{https://doi.org/10.1006/jabr.2000.8409}.
\bibitem{KMS22}V. Koppen, E. Meir, C. Schweigert, On isotypic decompositions for non-semisimple Hopf algebras, Algebr. Represent. Theory 25 (2) (2022) 447--475. \url{https://doi.org/10.1007/s10468-021-10029-x}.
\bibitem{LR88}R. G. Larson, D. E. Radford, Finite dimensional cosemisimple Hopf algebras in characteristic 0 are semisimple, J. Algebra 117 (2) (1988) 267--289. \url{https://doi.org/10.1016/0021-8693(88)90107-X}.
\bibitem{L97}M. Lorenz, Representations of finite-dimensional Hopf algebras, J. Algebra 188 (2) (1997) 476--505. \url{https://doi.org/10.1006/jabr.1996.6827}.
\bibitem{L98}M. Lorenz, On the class equation for Hopf algebras, Proc. Amer. Math. Soc. 126 (10) (1998) 2841--2844. \url{https://doi.org/10.1090/S0002-9939-98-04392-5}.
\bibitem{Mas95Two}A. Masuoka, Semisimple Hopf algebras of dimension $2p$, Comm. Algebra 23 (5) (1995) 1931--1940. \url{https://doi.org/10.1080/00927879508825319}.
\bibitem{Mas95Eight}A. Masuoka, Semisimple Hopf algebras of dimension 6, 8, Israel J. Math. 92 (1--3) (1995) 361--373. \url{https://doi.org/10.1007/BF02762089}.
\bibitem{Mas95Cube}A. Masuoka, Self-dual Hopf algebras of dimension $p^3$ obtained by extension, J. Algebra 178 (3) (1995) 791--806. \url{https://doi.org/10.1006/jabr.1995.1378}.
\bibitem{Mas96}A. Masuoka, The $p^n$ theorem for semisimple Hopf algebras, Proc. Amer. Math. Soc. 124 (3) (1996) 735--737. \url{https://doi.org/10.1090/S0002-9939-96-03147-4}.
\bibitem{MonS}S. Montgomery, Hopf algebras and their actions on rings, CBMS Regional Conference Series in Mathematics, vol. 82, American Mathematical
Society, Providence, RI, 1993. \url{https://doi.org/10.1090/cbms/082}.
\bibitem{Nat99}S. Natale, On semisimple Hopf algebras of dimension $pq^2$, J. Algebra 221 (1) (1999) 242--278. \url{https://doi.org/10.1006/jabr.1999.7968}.
\bibitem{Nat07Memo}S. Natale, Semisolvability of semisimple Hopf algebras of low dimension, Mem. Amer. Math. Soc. 186 (874) (2007), viii+123 pp. \url{https://doi.org/10.1090/memo/0874}.
\bibitem{Ng02}S.-H. Ng, Non-semisimple Hopf algebras of dimension $p^2$, J. Algebra 255 (1) (2002) 182--197. \url{https://doi.org/10.1016/S0021-8693(02)00139-4}.
\bibitem{NgPQ04}S.-H. Ng, Hopf algebras of dimension $pq$, J. Algebra 276 (1) (2004) 399--406. \url{https://doi.org/10.1016/j.jalgebra.2003.11.008}.
\bibitem{Ng2p05}S.-H. Ng, Hopf algebras of dimension $2p$, Proc. Amer. Math. Soc. 133 (8) (2005) 2237--2242. \url{https://doi.org/10.1090/S0002-9939-05-07804-4}.
\bibitem{NgPQII}S.-H. Ng, Hopf algebras of dimension $pq$, II, J. Algebra 319 (7) (2008) 2772--2788. \url{https://doi.org/10.1016/j.jalgebra.2007.08.003}.
\bibitem{NW26}S.-H. Ng, X. Wang, On stable twist $2$-groups of quasi-Hopf algebras, arXiv:2610.40228 [math.QA], 2026. \url{https://arxiv.org/abs/2610.40228}.
\bibitem{NZ89}W. D. Nichols, M. B. Zoeller, A Hopf algebra freeness theorem, Amer. J. Math. 111 (2) (1989) 381--385. \url{https://doi.org/10.2307/2374514}.
\bibitem{Stefan99}D. \c{S}tefan, Hopf algebras of low dimension, J. Algebra 211 (1) (1999) 343--361. \url{https://doi.org/10.1006/jabr.1998.7602}.
\bibitem{SVO98}D. \c{S}tefan, F. Van Oystaeyen, Hochschild cohomology and the coradical filtration of pointed coalgebras: applications, J. Algebra 210 (2) (1998) 535--556. \url{https://doi.org/10.1006/jabr.1998.7594}.
\bibitem{Wil88}R. E. Williams, Finite dimensional Hopf algebras, Ph.D. thesis, Florida State University, Tallahassee, FL, 1988.
\bibitem{Zhu94}Y. Zhu, Hopf algebras of prime dimension, Int. Math. Res. Not. 1994 (1) (1994) 53--59. \url{https://doi.org/10.1155/S1073792894000073}.
\end{thebibliography}
\end{document}